\documentclass[mathpazo]{cicp}
\usepackage{multirow}
\usepackage{graphicx}
\usepackage{pgf}

\newcommand{\be}{\begin{equation}}
\newcommand{\ee}{\end{equation}}

\newcommand{\bM}{\begin{displaymath}}
\newcommand{\eM}{\end{displaymath}}

\newcommand{\rsup}[1]{\mbox{\tiny #1}}

\newcommand{\bmath}{\begin{displaymath}}
\newcommand{\emath}{\end{displaymath}}

\newcommand{\bfg}{\begin{figure}}
\newcommand{\efg}{\end{figure}}

\newcommand{\bc}{\begin{center}}
\newcommand{\ec}{\end{center}}

\newcommand{\vect}[1]{{\bf #1}}

\newcommand{\al}{\alpha}
\newcommand{\del}{\delta}
\newcommand{\eps}{\epsilon}
\newcommand{\lam}{\lambda}
\newcommand{\pa}{\partial}

\def\B{{\mathcal B}}
\def\N{{\mathcal N}}
\def\T{{\mathcal T}}

\def\laplace{{\triangle}}

\begin{document}

\title{A Simple Method for Computing Singular or Nearly Singular
Integrals on Closed Surfaces}

\author[J. T. Beale et al.]{J. Thomas Beale\affil{1}\comma\corrauth,
                            Wenjun Ying\affil{2}, and
                            Jason R. Wilson\affil{3}   }
\address{\affilnum{1}\
Department of Mathematics, Duke University, Durham, NC 27708-0320, 
USA.\\
    \affilnum{2}\
Department of Mathematics, MOE-LSC and Institute of Natural Sciences,\\
\hspace*{2em} Shanghai Jiao Tong University, Minhang, Shanghai 200240, P. R. China.\\
    \affilnum{3}\ Mathematics Department, Virginia Tech, Blacksburg, VA 24061-0123, USA.
} 
\emails{{\tt beale@math.duke.edu} (J.~T.~Beale),
        {\tt wying@sjtu.edu.cn} (W.~Ying),
        {\tt jasonwil@math.vt.edu} (J.~R.~Wilson)  }

\begin{abstract}
We present a simple, accurate method for computing singular or nearly
singular integrals on a smooth, closed surface, such as layer potentials
for harmonic functions evaluated at points on or near the surface.
The integral is computed with a regularized kernel and corrections
are added for regularization and discretization, which are found
from analysis near the singular point. The surface integrals are
computed from a new quadrature rule using surface points which project
onto grid points in coordinate planes.  The method does not require
coordinate charts on the surface or special treatment of the
singularity other than the corrections.  The accuracy is about
$O(h^3)$, where $h$ is the spacing in the background grid,
uniformly with respect to the point of evaluation, on or near the surface.
Improved accuracy is obtained for points on the surface.
The treecode of Duan and Krasny for Ewald summation is used to perform sums.
Numerical examples are presented with a variety of surfaces.
\end{abstract}

\ams{65R20, 65D30, 31B10, 35J08}

\keywords{
boundary integral method, layer potential, nearly singular integrals,
Laplace equation, surface integral, implicit surface
}

\maketitle

\section{\label{sec_intro}Introduction}

We present a simple, accurate method for computing singular or nearly
singular integrals defined on a smooth, closed surface
in three-space.  This method can be used to evaluate single
or double layer potentials for harmonic functions, or the
velocity and pressure in Stokes flow due to forces on a surface.
The point of evaluation could be on or near the surface.
To evaluate the integral,
the kernel is first replaced by a regularized
version.  A preliminary value is found using a new quadrature
rule for surface integrals which has the advantage that it
does not require coordinate systems or a triangulation
on the surface.  Instead we sum values of the integrand
over quadrature points which project onto grid points in coordinate
planes, in a way that would be high order accurate if
the integrand were smooth.  Corrections for
the regularization and discretization are then added to achieve
higher accuracy.  These corrections are given
by explicit formulas derived using asymptotic analysis near the
singularity as in \cite{Beale:2004}.
The resulting value of the
integral has $O(h^p)$ accuracy for $p < 3$, uniformly for points
of evaluation near the surface,
where $h$ is the grid spacing in $\mathbb R^3$.  For points on the
surface, the accuracy is significantly improved by using a special
regularization; see Section 3.3.  For efficient summation we use the
treecode algorithm of Duan and Krasny \cite{Duan:2000} designed for
kernels with Gaussian regularization.
The method presented here could be used,
for example, to find values of the potential at grid points in
$\mathbb R^3$ close to the surface.  It
should be applicable to computations with moving
surfaces for which good accuracy is needed without extensive
work to represent the surface at each time step.  Other kernels
could be treated by the same approach, and the more accurate
version of the method for
computing values on the surface could be
used for a variety of problems which can be formulated as
integral equations.

The present approach is an improvement and extension of the
grid-based boundary integral method
of \cite{Beale:2004}.  In the earlier work the integral was replaced
by sums in coordinate charts using a partition of unity.  The
need for explicit coordinate systems requires knowledge of
the surface that might be difficult to obtain for a moving
surface.  Furthermore, if the coordinate system is too distorted, the
accuracy will be poor because the discretization error will fail to be
controlled by the regularization.   Here we avoid these disadvantages
by using a more direct rule for computing surface integrals,
which was introduced
for smooth integrands in \cite{Wilson:2010}.  This quadrature rule
uses projections on coordinate planes rather than coordinate charts.
Given a rectangular grid in $\mathbb R^3$,
quadrature points are chosen as those points on the surface which
project onto grid points in the coordinate planes, for which the
normal to the surface has direction away from the plane.  Weights
for the quadrature points are found from a partition of unity on
the unit sphere, applied to the normal vector at the point.
The weight functions on the sphere are chosen universally, and do
not depend on the particular surface.  The resulting quadrature
rule for surface integrals has high order accuracy,
as allowed by the smoothness of the integrand and the surface.
In effect, the method uses the
existence of coordinate patches without having to refer to them
explicitly.  The quadrature points can be found efficiently if,
for example, the surface is given analytically or numerically
as the level set of a function.  Examples in \cite{Wilson:2010}
with smooth integrands illustrate the accuracy of this method 
with a variety of surfaces, including ones of large genus.

For a variety of problems in partial differential
equations, solutions can be written as integrals
over surfaces, using a known fundamental solution.  These
include harmonic functions, electromagnetic waves, and viscous
fluid flow modeled by the Stokes equations.  The specific representation
makes this formulation an attractive approach for numerical methods.
Often an integral equation on the surface must be solved, and much
attention has focused on such problems.  For the most familiar case,
an integral equation for the Dirichlet problem for the Laplacian,
using the double layer potential, it was proved in \cite{Beale:2004}
that the discrete version of the integral equation using the present method
has a unique solution
and that it converges to the exact solution under grid refinement.
The evaluation of the integral at points near
the surface, in contrast to evaluation on the surface,
generally requires extra care and can also be important in applications.
If values of a potential are needed at grid points in the computational domain, we can first
compute values at points near the surface directly as nearly singular integrals.  Once this is
done, values at other grid points can be found in a cheap way by inverting a discrete Laplacian,
as in \cite{Mayo85}.  In this paper the values obtained at the irregular grid points near the surface
using integration have
accuracy about $O(h^3)$, while those found at the regular points are about $O(h^2)$ accurate;
more accurate values could be found with more work.

The most widely used numerical technique for integral formulations
is the boundary element method (e.g. \cite{Atkinson:1997,
Hackbusch:1995, Pozrikidis:2002, SauterSchwab});
the boundary is triangulated and matrix elements are computed
for the integral operator on the surface using special quadrature
rules.  This method is especially useful for electromagnetic
problems in which the surface does not evolve and may have corners
and edges.  Direct quadrature based on triangulations is also
used, and can be accurate in the (exactly) singular case
for smooth surfaces.  Such methods have long been used for
modeling in chemistry and biology \cite{Zinchenko:1999, 
Pozrikidis:2001, Pozrikidis:2002}.  Some methods use corrected
quadrature weights \cite{Tornberg} or analytical evaluation of
the singular part, treating the remainder in a standard way
\cite{Helsing:2013, Wilson:2010}.
A careful direct quadrature (or Nystr\"om)
method, introduced in \cite{Bruno:2001} for electromagnetics,
and further developed in \cite{Biros:2006},
uses a partition of unity to reduce the integral to coordinate patches;
a special patch in polar coordinates is used near the singularity.
A different approach \cite{Graham:2002} is a spectral method,
using spherical harmonics, assuming the surface can be mapped to
a sphere.  This method, applied in \cite{Ganesh:2004, Biros:2011}, 
appears to be advantageous when there are many boundaries that are
not greatly deformed. For use with many vesicles, the authors
of \cite{ Biros:2011} report difficulties in keeping the method
of \cite{Bruno:2001, Biros:2006} accurate, perhaps because of the
cut-off function needed for the special patch near the singularity.
Careful integral methods have the promise of calculating integrals
very accurately even for complicated surfaces lacking smoothness;
see \cite{Helsing:2008, Bremer:2012, Klockner:2013, Hughes}.  An alternative
approach, the kernel-free boundary integral method
\cite{YingWang13,YingWang14}, replaces the
calculation of the integral by the solution of an interface problem
on a regular grid.  This approach has the advantage that explicit knowledge
of the integral kernel is not needed, and thus it can be applied to more
general partial differential equations.

One advantage of the present method is its simplicity.
Detailed information about the surface is not required.
No special treatment is needed near the singularity,
except for the corrections which are added after the summation,
using analytical formulas.
By design, the errors are uniform
with respect to the location; the additional work
needed for points close to the boundary is small compared
to that for points on the boundary or far away.
For an integral with sources on several boundaries,
the evaluation at a point on one surface might lead to
a nearly singular case, since it could require integration
over another surface which is close to the first.
Results in \cite{Ying:2013} showed that the
two-dimensional version \cite{Beale:2001} of the present
method works well in such cases.

In this paper we treat single and double layer potentials
for harmonic functions.  The approach can be extended to
integrals for Stokes flow as in \cite{Tlupova:2013};
see also \cite{Cortez:2001, Cortez:2005, Nguyen}.
This integration method could be
applied to problems with moving interfaces.
A discretization of the evolving surface must be chosen
and updated.  This is often
done with representative marker points and triangulation.
An alternative might be to use the level set method 
\cite{Osher:2003, Sethian:1998},
representing the surface as the zero set of a function
whose values are transported by a computed velocity.
The velocity would be needed at nearby grid points, in order
to update the level set function, rather
than on the surface, and the present integration method
is designed to be suitable for this purpose.

In Section \ref{sec_int_method} the 
quadrature rule for surface integrals is explained.
In Section \ref{sec_layer_potentials} the formulas for
calculating single and  double layer potentials are given,
including the simplified version for points on the surface.
A brief discussion of the error estimates is included.
Numerical examples illustrating the method with a variety
of surfaces are presented in Section \ref{sec_results}. 
Finally, some possible improvements
for this method are discussed briefly in Section 
\ref{sec_discussion}.
The code that produced the examples is available on request from
the first two authors.

\section{\label{sec_int_method}The quadrature method for surface integrals}

In this section we describe the computational method for integrals on
implicitly defined closed surfaces in three space dimensions.  The method is 
also applicable to closed curves in $\mathbb R^2$ and to hypersurfaces in
higher dimensions. Detailed proofs and extensive examples are given
in Wilson \cite{Wilson:2010}.
Our purpose is to evaluate the surface integral
\be
\label{surf-int}
I \,=\, \int_{\Gamma} f(\vect y) \, d S_{\vect y}
\ee 
We assume $\Gamma$ is a $C^{2m+1}$ surface,  $m \geq 1$,
and, for now, $f \in C^{2m}(\Gamma)$.
The method exploits the spectral convergence of the trapezoidal rule 
without the need for the user to generate a set of overlapping coordinate
patches and associated partition of unity. Rather than covering $\Gamma$
with overlapping rectangular patches, we cover $\Gamma$ with certain
overlapping surface sets. We define the subsets
\bM
\Gamma_i = \{\, \vect x \in \Gamma \, : \, |\vect n (\vect x) \cdot \vect e_i | > 0 \, \},
\quad i = 1,2,3
\eM 
where $\vect n$ is the unit outward normal at $\vect x$ and 
$\{\vect e_i\}_{i=1}^3$ is the standard basis for $\mathbb R^3$.
Thus $\Gamma_i$ contains all points in $\Gamma$ where $\vect n$ is
not orthogonal to $\vect e_i$.

The integration method uses a high order,
patch--independent quadrature formula for integrals
with integrands that vanish outside of a compact subset of one $\Gamma_i$. 
To handle general integrands $f$, we introduce a partition of unity to find
functions $\{f^i\}_{i=1}^3$ such that $f(\vect x) = \sum_{i=1}^3 f^i(\vect x)$
for all $\vect x \in \Gamma$ and $f^i$ vanishes outside a compact subset of 
$\Gamma_i$. We first design a universal partition of unity on the unit
sphere, $S = \{ \, \vect u \in \mathbb R^3 \, : \, |\vect u| = 1 \, \}$. We 
start with the smooth bump function defined as $b(r) = e^{r^2/(r^2 - 1)}$ 
for $|r| < 1$ and $b(r) = 0$ otherwise. Next we choose a fixed angle 
$\theta$ with  $\cos^{-1}(1/\sqrt{3}) < \theta < \pi/2$. For each $i \in \{1, 2, 3\}$
and $\vect u \in S$ we define
\bM
w^i(\vect u) = \cos^{-1}(|\vect u \cdot \vect e_i|)\,, \qquad
\sigma^{i,\theta}(\vect u) = 
\frac{b(w^i(\vect u)/\theta)}{\sum\limits_{j=1}^3 b(w^j(\vect u)/\theta)}
\eM 
Because $\theta > \cos^{-1}(1/\sqrt{3})$, the sum is always positive.
Furthermore
\begin{enumerate}
\item
For each $i = 1, 2, 3$, we have $\sigma^{i,\theta} \in C^{\infty}(S)$; 
\item
For all $\vect u \in S$, we have $\sum_{i=1}^3 \sigma^{i,\theta}(\vect u) = 1$;
\item
For each $i = 1, 2, 3$, the function $\sigma^{i,\theta}$ vanishes outside the
compact subset 
\bM
S_{i,\theta} = \{\, \vect u \in S \, : \, 
|\vect u \cdot \vect e_i | \geq \cos{\theta} \, \}
\eM 
\end{enumerate}
To make use of the above partition of unity on the sphere
for a general surface $\Gamma$,
we apply it to the unit normal $\vect n (\vect x) \in C^{2m}(\Gamma)$.
The composition functions $\zeta^{i,\theta} = \sigma^{i,\theta} \circ \vect n$
on $\Gamma$ satisfy
\begin{enumerate}
\item
For each $i = 1, 2, 3$, we have $\zeta^{i,\theta}  \in C^{2m}
(\Gamma)$; 
\item
For all $\vect x \in \Gamma$, we have $\sum_{i=1}^3 \zeta^{i,\theta}(\vect x)  = 1$;
\item
For each $i = 1, 2, 3$, the function $\zeta^{i,\theta}  $ vanishes 
outside the compact subset of $\Gamma_i$ given by 
\be
\label{Gammaitheta}
\Gamma_{i,\theta} = \{\, x \in \Gamma \, : 
\, |\vect n (\vect x) \cdot \vect e_i | \geq \cos{\theta} \, \}
\ee 
\end{enumerate}

Using the partition of unity, we obtain the exact formula
\be
\int_{\Gamma} f(\vect y) \, d S_{\vect y} 
= \sum\limits_{i=1}^3 \int_{\Gamma_i} \zeta^{i,\theta}(\vect y)
	f(\vect y) \, d S_{\vect y} 
\ee 
where the integrand $f^i(\vect y) = \zeta^{i,\theta}(\vect y) f(\vect y)$
vanishes outside the compact subset $\Gamma_{i,\theta}$ of $\Gamma_i$. 
Finally, the surface integral (\ref{surf-int}) can be approximated by the
numerical quadrature 
\be
\label{smoothquad}
I_h \,=\, h^2 \sum\limits_{i=1}^3 \sum\limits_{\vect x \in R_{h,i,\theta}}
\frac{\zeta^{i,\theta}(\vect x) f(\vect x)}
{|\vect n (\vect x) \cdot \vect e_i|}  
\,=\, h^2 \sum\limits_{i=1}^3 \sum\limits_{\vect x \in R_{h,i,\theta}}
\frac{\sigma^{i,\theta}(\vect n (\vect x)) f(\vect x)}
{|\vect n (\vect x) \cdot \vect e_i|}  
\ee 
where
\bM
R_{h, i, \theta} = \{ \, \vect x \in \Gamma \, : \, |\vect n (\vect x) \cdot 
\vect e_i | \geq \cos{\theta} \, 
\; \mbox{and} \;  p^i(\vect x) \in h \, \mathbb Z^2 \, \} 
\eM 
and $p^i: \mathbb R^3 \mapsto \mathbb R^2$ is 
the projection function defined by 
$ p^1(\vect x) = (x_2, x_3)$,
$ p^2(\vect x) = (x_1, x_3)$,
$ p^3(\vect x) = (x_1, x_2)$.
$R_{h, i, \theta}$ consists of those points in $\Gamma_{i,\theta}$ that
project to grid points in the corresponding plane.
The weights $1/ |\vect n (\vect x) \cdot \vect e_i|$ correspond
to the area elements of the inverse projections.

It is proved in \cite{Wilson:2010} that $I_h - I = O(h^{2m})$, i.e., the quadrature rule
(\ref{smoothquad}) is high order accurate, provided
the surface $\Gamma$ is $C^{2m+1}$ and $f$ is $C^{2m}$, $m \geq 1$;
see Lemma 1, Theorem 2, pp. 9--10, and Lemma 9, Theorem 10, pp. 25--27.
In effect the trapezoidal rule applies on coordinate patches covering each
$\Gamma_{i,\theta}$.
The spacing $h$ must be small enough to resolve the surface.
Assume $\Gamma$ is defined as the set $\phi = 0$ for some function
$\phi$ on an open subset of $\mathbb R^3$. 
If $C_1 = \min {|\nabla\phi|}$ and
$C_2 = \max {\|D^2\phi\|}$ near the surface, we need
$h < h_0 \equiv 2C_1\cos{\theta}/C_2$.  Thus if the curvature is large, $h$ must be small.  The method works if $\phi$ is known
only at grid points; see p. 30 of \cite{Wilson:2010}.

The use of this rule requires finding the points in $R_{h,i,\theta}$.  Given a grid
point $\hat{x}$ in a coordinate plane, there may be several points in $\Gamma$ which
project to $\hat{x}$, but they are well separated because of the normal condition
in (\ref{Gammaitheta}).
Consequently, as shown in \cite{Wilson:2010}, pp. 11--18,
a simple line search algorithm can be used to locate the quadrature
points if $h < h_0$:  Briefly, to find
points in $R_{h,3,\theta}$, for each $(j_1,j_2)$ and $j_3$, check whether $\phi$ has a root ${\bf x} = (j_1h,j_2h,x_3)$ with $j_3h \leq x_3 \leq (j_3+1)h$.  If so, find a root.
If the root is in $R_{h,3,\theta}$, it is unique.  If it is not in $R_{h,3,\theta}$,
reject it.  In either case go to the next $j_3$. 
 The validity of this algorithm for $h < h_0$ is
proved in \cite{Wilson:2010}.  

\section{\label{sec_layer_potentials}Evaluation of the layer potentials}

We describe the procedure for computing a single or double layer potential at an
arbitrary point, the most difficult case being a location off the surface but close by.
For the case of a point on the surface we 
give more special versions with improved accuracy.
Finally we discuss error estimates.

\subsection{The single layer potential} 

The single layer potential on $\Gamma$ determined by a density function
$\psi$ is 
\be
\label{single-layer-v}
v(\vect x) = \int_{ \Gamma} G(\vect y - \vect x) \psi(\vect y)
	\, d S_{\vect y}
\ee 
where $G$ is the fundamental solution for the Laplacian, $G(x) = -1/4\pi |\vect x|$.
We suppose $\Gamma$ is the boundary of a bounded domain $\Omega$.
To evaluate $v$ for $\vect x$  close to $\Gamma$,
we replace $G$ with a smoothed, or regularized, version
\be
\label{Gdelta}
G_{\delta}(\vect y) = G(\vect y) \, \text{erf}(|\vect y|/\delta)
  = - \frac{\text{erf}(|\vect y|/\delta)}{4\pi |\vect y|} 
\ee 
where $\text{erf}$ is the error function and
$G_{\delta}(\vect 0) = -\pi^{-3/2}(2\delta)^{-1}$. 
The resulting error in the integral is $O(\delta)$. 
Typically $1 \leq \delta/h \leq 2$.
We first compute the regularized integral
\be
\label{reg-single-layer}
v_{\delta}(\vect x) =  \int_{ \Gamma} G_{\delta}(\vect y - \vect x) 
\psi(\vect y) \, d S_{\vect y}.
\ee 
using the method of Sec. \ref{sec_int_method}.
The value obtained
is not close to $v(\vect x)$ because of the near singularity.
We add corrections for the regularization and discretization
to improve the accuracy.  Other regularizations could be used, but 
(\ref{Gdelta}) has the advantages that it is simple, $G_\delta - G$ decays
rapidly in the far field, and manageable formulas can be found for the
corrections described below.  

To obtain the corrections we first find
$\vect z \in \Gamma$, the closest point on the surface to
$\vect x$, and set $\vect x = \vect z + b \, \vect n$. 
Here $\vect n$ is the unit outward normal to $\Gamma$ at $\vect z$; $b < 0$ if 
$\vect x \in \Omega$ and $b > 0$ if $\vect x \in \Omega^c$ (the complement
of the closure of $\Omega$ in $\mathbb R^3$). 
The correction for regularization of the single layer potential is
\be
\label{T1} 
\T_1 =  \frac{\delta}{2} \bigl (1 + H \lambda \, \delta \bigr ) \, 
\psi(\vect z) \, 
  \biggl [
  |\lambda| \, \text{erfc} |\lambda| - \frac{e^{-\lambda^2}}{\sqrt{\pi}}
  \biggr ] .
\ee 
Here $\text{erfc}(r) = 1 - \text{erf}(r)$,
$\lambda = b/\delta$ and $H$ is the mean curvature at $\vect z$, 
$H = (\kappa_1 + \kappa_2)/2$, where $\kappa_1$ and $\kappa_2$ are the 
principal curvatures.  Formulas for computing needed geometric quantities
such as $H$ are given in Appendix B, and the sign convention for $H$ is
explained.

The discretization correction is a rapidly convergent infinite sum resulting
from the Poisson summation formula, applied in each $\Gamma_i$, $i = 1,2,3$.  Let
\be
E(p,q) = e^{2pq}\text{erfc}(p+q) + e^{-2pq}\text{erfc}(-p+q) \,.
\ee
Also let $Q = \{n = (n_1,n_2) \in \mathbb{Z}^2 \,: n_2 > 0 \;\; \mbox{or} \;\;
(n_2 = 0\; \mbox{and} \; n_1 > 0)\}$.  Now suppose $\vect z$ lies in a system of
coordinates, say $\al = (\al_1,\al_2)$; in our case,  $\vect z \in \Gamma_{k,\theta}$
for one or more of $k=1,2,3$,
and the coordinates at $\vect z$ are $p^k(\vect z) \in \mathbb R^2$.
Let $g^{ij}$ be the inverse metric
tensor at $\vect z$, and for $n \in Q$, define $\|n\|^2 = \sum_{i,j} g^{ij}n_in_j$.
Also write the coordinates of $\vect z$ as $(\al_1,\al_2) = (m_1,m_2)h + (\nu_1,\nu_2)h$
where $m_1, m_2$ are integers and $0 \leq \nu_1,\nu_2 < 1$.
Here $(\nu_1,\nu_2) = \nu^{(k)}$ and $\|n\|$  depend on the
choice of $k$.  The discretization correction is
\be
\label{T2}
\T_2 = \frac{h}{4\pi} \psi(\vect z) \sum_{k=1}^3 \sum_{n \in Q} \zeta^{k,\theta}(\vect z)
            \cos(2\pi n\cdot\nu^{(k)}) \frac{1}{\|n\|} E(\lambda, \pi \delta \|n\|/h)
\ee
with $\zeta^{k,\theta}$ as in Sec. \ref{sec_int_method}.
Finally, the computed value $v$ of the integral (\ref{single-layer-v}) is
$v \approx v_\delta + \T_1 + \T_2$, where $v_\delta$ is the
value of (\ref{reg-single-layer}) obtained by the quadrature rule.

\subsection{The double layer potential} 

The double layer potential has the form
\be
\label{double-layer-v}
w(\vect x) = \int_{ \Gamma} 
\frac{ \partial G(\vect y - \vect x)}{\partial \vect n_{\vect y}}
\varphi(\vect y)
\, d S_{\vect y}
\ee 
It is discontinuous at $\Gamma$.  If $\vect x$
is close to $\Gamma$, we find the closest point $\vect z$ and distance $b$ as before.  
We use Green's identities to reduce the singularity and then regularize the kernel,
obtaining
\be
\label{reg-double-layer}
w_{\delta}(\vect x) = 
\int_{ \Gamma} 
\frac{\partial G_{\delta}(\vect y - \vect x)}{\partial \vect n_{\vect y}} 
\bigl [ \varphi(\vect y) - \varphi(\vect z) \bigr ] \, d S_{\vect y}
+ \chi \varphi(\vect z)
\ee 
Here $\chi = 1$ for $x \in \Omega$, $\chi = 0$ on $\Omega^c$,
$\chi = \textstyle{\frac12}$ on $\Gamma$.
To form $\pa G_\delta/\pa \vect n$ we use the gradient of the smooth
function $G_\delta$ introduced in (\ref{Gdelta}),
\be
\label{gradGdel}
\nabla G_{\delta}(\vect y) = \nabla G(\vect y)s(|\vect y|/\delta)
    =   \frac{\vect y}{4 \pi |\vect y|^3}s(\vect y|/\delta)
\ee
with 
\be
\label{sforgrad}
s(r) = \text{erf}(r) - \frac{2}{\sqrt{\pi}} \, r e^{-r^2}
\ee 
We compute the integral in (\ref{reg-double-layer}) as in Sec. \ref{sec_int_method}
and again add corrections.
The regularization correction for
(\ref{reg-double-layer}) is
\be
\label{N1} 
\N_1 = \delta^2 \, (\laplace_S \varphi) \, \frac{\lambda}{4} 
\biggl [
|\lambda| \, \text{erfc}|\lambda| - \frac{e^{-\lambda^2}}{\sqrt{\pi}}
\biggr ]  .
\ee 
Here
$\laplace_S \varphi$ is the surface Laplacian of $\varphi$ at $\vect z$, which is 
expressed in coordinates as
\bM
\laplace_S \varphi = \frac{1}{\sqrt{g}} \sum\limits_{i,j =1}^2 
\frac{\partial}{\partial \al_j} 
\biggl (
    \sqrt{g} \, g^{ij} \frac{\partial \varphi}{\partial \al_i} 
    \biggr ).
\eM
with $g = \det{g_{ij}}$.
The discretization correction is similar to $\T_2$ but involves 
$\pa\varphi/\pa \al_r$, $r = 1,2$, the coordinate derivatives of $\varphi$
evaluated at $\vect z$.  In our case, on $\Gamma_{k,\theta}$, $(\al_1,\al_2)$ are
the two components of $\vect x = (x_1,x_2,x_3)$ other than $x_k$, and
we write these derivatives as $\pa_r^{(k)}\varphi(\vect z)$.
The correction is
\be
\label{N2}
 \N_2 = -\frac{\delta\lambda}{2} \sum_{k=1}^3 \sum_{r=1}^2 c_r^{(k)} 
            \zeta^{k,\theta}(\vect z) \pa_r^{(k)}\varphi(\vect z)
\ee
where
\be 
\label{N2c}
 c_r^{(k)} = \sum_{n \in Q} \sum_{s=1}^2 \sin(2\pi n\cdot\nu^{(k)}) 
    \frac{g^{rs}n_s}{\|n\|} E(\lambda, \pi \delta \|n\|/h)
\ee
with $\|n\|$ as before.
The computed value of (\ref{double-layer-v}) is $w \approx w_\delta + \N_1 + \N_2$.

\subsection{The potentials evaluated on the surface}

We now treat the important special case of 
evaluation at a point $\vect x$ on the surface $\Gamma$.  
For this case, in contrast to the nearly singular case, it is not difficult
to modify the regularized kernels to have higher accuracy by imposing moment conditions.  Thus no corrections are needed for regularization.  The method
is easier to use than in the general case, and the error is typically smaller,
as seen in the examples in Section 4.
A strategy for producing these improved kernels from the ones
already chosen in (\ref{Gdelta}) and (\ref{gradGdel}) is described in \cite{Bealeww} for the single layer
and in \cite{Beale:2004} for the double layer.  For $G_\delta$ as defined in
(\ref{Gdelta}), the error in regularization, i.e. the difference between the integrals in (\ref{reg-single-layer}) and (\ref{single-layer-v}), is $O(\delta)$.
For the new version of $G_\delta$ the error is $O(\delta^5)$ for the special
case of evaluation on the surface, and similarly for the double layer. 

To evaluate the single layer potential (\ref{single-layer-v})
at $x \in\Gamma$, we use the new version of $G_\delta$ with $O(\delta^5)$
accuracy,
\be
\label{Gdelta5}
G_{\delta}(\vect y) = - \frac{s(|\vect y|/\delta)}{4\pi |\vect y|} \,,\quad
        s(r) = \text{erf}(r) + \frac{2}{3\sqrt{\pi}}(5r - 2r^3)\,e^{-r^2} 
\ee 
and $G_{\delta}(\vect 0) = -(4/3)\pi^{-3/2}\delta^{-1}$. In place of the
corrections (\ref{T1}, \ref{T2}) we have $\T_1 = 0$ and
\be
\label{T25}
\T_2 = \frac{\delta}{\pi} \psi(\vect z) \sum_{k=1}^3 \sum_{n \in Q} \zeta^{k,\theta}(\vect z)
            \cos(2\pi n\cdot\nu^{(k)})\,F(\xi)
\ee
where, with $\xi = 2\pi\|n\|\delta/h$ and  $\|n\|$ as before,
\be
F(\xi) = \frac{\pi}{\xi} \text{erfc}(\xi/2)
      + \pi^{1/2}\frac{\delta}{h}\,e^{-\xi^2/4}\,\left(1 + \frac{\xi^2}{6}\right)
\ee
The derivation of (\ref{Gdelta5}) is similar to that in Sec. 2 of \cite{Bealeww} for an
$O(\delta^3)$ version, with the
formula (\ref{T25}) corresponding to (3.28), (2.23) in \cite{Beale:2004}.

To evaluate the double layer potential (\ref{double-layer-v}) at $\vect x \in\Gamma$, we use
(\ref{reg-double-layer}) with $\chi = \frac12$ and 
$\nabla G_\delta$ of the form (\ref{gradGdel}) but with (\ref{sforgrad}) replaced by
\be
\label{sforgrad5}
s(r) = \text{erf}(r) - \frac{2}{\sqrt{\pi}} \, \left(r - \frac{2r^3}{3}\right) e^{-r^2}.
\ee 
In this case no corrections are needed, that is, $\N_1 = \N_2 = 0$.  Formula (\ref{sforgrad5})
was derived in \cite{Beale:2004}, p. 607.

\subsection{Error analysis} 

For the general case of points close to $\Gamma$, the
corrections $\N_1$, $N_2$ for the double layer were derived in \cite{Beale:2004}.
Those for the single layer can be found similarly; we include a brief derivation
of $\T_1$ in Appendix A.
The discretization corrections are based on the Poisson Summation Formula.
After applying both corrections to either the single or double layer potential,
the remaining error has the form (cf. Theorem 1.2 in \cite{Beale:2004}) 
\be 
\label{errorestimate}
\epsilon \,\leq\, C_1\delta^3 + C_2 h^2 e^{-c_0 (\delta/h)^2}
\ee
as $\delta$, $h \to 0$, assuming $\delta/h$ is bounded below, with $C_1, C_2$ depending on derivatives of the surface and density functions.  The two terms
represent the regularization error and discretization error, respectively.
The constant $c_0$ is determined by the
choice of local coordinate system.  It is important for accuracy that $c_0$ does not become small, so that the sums in the discretization corrections converge rapidly, and so that the second error term in (\ref{errorestimate})
is comparable to the first in practice.
It was shown in \cite{Beale:2004}, Sec. 3, that the estimate (\ref{errorestimate}) holds provided
$$ c_0 \,<\, \pi^2 \gamma^2 \,, \quad 
  \gamma^2 \,\equiv\, \min_{|k| = 1} \sum_{i,j} g^{ij} k_i k_j \,. $$
(A factor $1/2$ in $c_0$ in \cite{Beale:2004} was arbitrary.)
Here $g^{ij} = (g_{ij})^{-1}$ is the inverse metric tensor.  With coordinates
$(\alpha_1,\alpha_2)$, $g_{ij} = T_i\cdot T_j$,
where $T_j = \pa{\bf x}/ \pa \alpha_j$, $j = 1,2$ are the
tangent vectors.
If the coordinate system distorts distances significantly, $\gamma$ could be small, and thus the accuracy of the method depends on the choice of coordinates.

In our case, the coordinate systems are those determined by the projections.
For $\Gamma_3$, the coordinates are $x_1,x_2$, with $x_3 = f(x_1,x_2)$,
and similarly for $\Gamma_1$, $\Gamma_2$.  From the expression for $g^{ij}$ in
Appendix B, it is not difficult to see that
in this case $\gamma = (1 + f_1^2 + f_2^2)^{-1/2} = 
|\vect n\cdot \vect e_3| \geq \cos{\theta}$ in $\Gamma_{3,\theta}$, where 
$\theta$ is the angle chosen in Sec. \ref{sec_int_method}.
Thus (\ref{errorestimate}) holds with $c_0 \approx \pi^2\cos^2{\theta}$,
and the exponential in (\ref{errorestimate}) can be made quite small.  For example, if
$\delta = 2h$ and $\theta = 60^o$,
the exponential is $.00005$; if $\delta = 2.5h$ and $\theta = 70^o$ it is $.0007$.
In practice the accuracy is about $O(h^3)$ for usual values of $h$,
with proper choice of parameters.  Alternatively, we could take $\delta = Ch^q$
for any $q < 1$ and thereby obtain convergence as $h \to 0$ with order $O(h^{3q})$.

For the sums in (\ref{T2}), (\ref{N2c}) we only need a few terms because of
the rapid decay as $n$ increases.  In the corrections we
may evaluate $H$ and $\laplace_S\varphi$ not at the closest point $\vect z$
but rather at a neighboring grid point, using formulas in Appendix B.
The $O(h)$ errors in these quantities do not change the order of accuracy of the corrections.
Similarly, in the discretization corrections,
the $g^{ij}$ and $\partial^{(k)}_r\varphi$
only need to be computed within $O(h)$.

For the special case with $\vect x \in \Gamma$, the first term in the error estimate
(\ref{errorestimate}) improves to $\delta^5$.  (See Theorem 1.1 of \cite{Beale:2004}
for the double layer.)  Thus we can take $\delta/h$ somewhat
larger to reduce the discretization error in the second term.  If
$\delta = Ch^q$, $q < 1$, then the error is $O(h^{5q})$ as $q \to 0$.

\section{\label{sec_results}Numerical Results}

This section presents examples evaluating
the sum of a double layer and single layer potential on five different surfaces,
\bM
u(\vect x) = w(\vect x) + v(\vect x) = 
\int_{\Gamma} \frac{\partial G(\vect y, \vect x)}{\partial \vect n_{\vect y}}
\varphi(\vect y) \, ds_{\vect y} 
- 
\int_{\Gamma} G(\vect y, \vect x) \psi(\vect y) \, ds_{\vect y} .
\eM 
In all the examples the potential $u$ is chosen to be
\bM
u(x, y, z) = 
\begin{cases}
(\sin x + \sin y) e^z & \text{if $(x, y, z) \in \Omega$}  \\
0 & \text{if $(x, y, z) \in \Omega^c$}
\end{cases}.
\eM 
The densities $\varphi$ and $\psi$ are determined by the jumps in $u$ and $\pa u/\pa n$.
The integrals are calculated given these densities, and the result is compared with
the exact $u$.  This choice of test problem allows us to have an exact solution with an
arbitrary surface.

In each example, the domain $\Omega$ is embedded into a cubic 
box $\B = (-L,L)^3$ with $L = 1.1$.
The box is partitioned into a uniform grid
$\T_h$ with mesh parameter $h = 2L/N$, the width of a grid cell.  
We call a grid node {\it irregular} if the stencil of the second-order
Laplacian $\Delta_h$ crosses the boundary $\Gamma$; otherwise it is
{\it regular}. 

The numerical values of $u(\vect x) = w(\vect x) + v(\vect x)$ are first computed
at the irregular grid nodes, as well as the neighboring nodes in their
stencil, using the procedure of Sec. \ref{sec_layer_potentials}.
The sums $w_\delta$, $v_\delta$
for the smoothed potentials are found and the corrections are added.
The summation is done using a slight modification of the treecode algorithm of
Duan and Krasny \cite{Duan:2000}, which was designed for
use with Ewald summation.  (The kernel in their code has a factor of erfc
rather than erf, and their solutions are periodic rather than in free space.)
In the treecode we chose the degree of Taylor polynomials
$p = 12$, the separation parameter $s = .5$, and the capacity, or maximum
number of points in a leaf, $N_0 = 20$. 

Having calculated $u(\vect x)$ at grid nodes $\vect x$ near $\Gamma$,
we can now find values at all the regular grid nodes of $\T_h$
by inverting the discrete Laplacian, using a procedure suggested in \cite{Mayo85}.
Let $u_{h,\delta}$ denote the value of $u$ already computed at the nodes close
to $\Gamma$
as nearly singular integrals.  We formulate a Poisson problem for
an approximation $u_h$ on $\T_h$ to the exact $u$,
\bM
\Delta_h u_h = \begin{cases}
\Delta_h u_{h, \delta} & \text{at irregular grid nodes} \\
0 & \text{at regular grid nodes} 
\end{cases}.
\eM
We solve for $u_h$ with a fast Poisson
solver on the box $\B$ with zero boundary condition.
At the regular nodes, the truncation error is $O(h^2)$, since the exact solution
is a smooth, harmonic function away from $\Gamma$.  For the irregular nodes,
$u_{h, \delta} - u$ is about $O(h^3)$, so that $\Delta_h u_{h, \delta} -
\Delta_h u = O(h^3/h^2) = O(h)$.  In summary, $\Delta_h (u_h - u)$ is $O(h^2)$
at regular nodes and $O(h)$ at irregular nodes.  Despite the first order truncation error near $\Gamma$, the resulting $u_h$ is second order accurate, i.e., 
$u_h - u = O(h^2)$, uniformly on $\T_h$, a fact proved in 
\cite{bealay}.  While this method is efficient, more accurate values
could be computed as integrals or otherwise.

We also computed the values of $u(\vect x)$ at grid nodes $\vect x$ on $\Gamma$
using the method in Section 3.3.  (The exact value is the average of the inside
and outside limits.)  In this case we used direct summation to provide
an unambiguous test of the accuracy.  

For Example 1, the surface $\Gamma$ is an ellipsoid given by
\bM \frac{x^2}{a^2} + \frac{y^2}{b^2} + \frac{z^2}{c^2} = 1 \eM
with $a = 1$, $b = .8$, $c = .6$, and rotated by an orthogonal matrix to test the
effect of grid alignment.
Example 2 is a thinner ellipsoid, with $a =1$, $b = c = .4$,
without rotation.
Example 3 is a torus 

\bM  (\sqrt{x^2 + y^2} - c)^2 + z^2 = a^2 \eM
with $a = .3$ and $c = .7$.
Example 4 is a molecular surface with four atoms, similar to one of the definitions in \cite{tms},
\bM  \sum_{k=1}^4 \exp(|\vect x - \vect x_k|^2/r^2) = c  \eM
Here the centers are
$(\sqrt{3}/3,0,-\sqrt{6}/12)$, $(-\sqrt{3}/6,.5,-\sqrt{6}/12)$,
$(-\sqrt{3}/6,-.5,-\sqrt{6}/12)$, $(0,0,\sqrt{6}/4)$
and $r = .5$, $c = .6$.
Example 5 is a surface obtained by revolving a Cassini oval,
\bM (x^2 + y^2 + z^2 + a^2)^2 - 4a^2(x^2 + y^2) = b^4  \eM
with $a = .65$ and $b = .7$. 

\begin{figure}
  \centering
 \includegraphics[width=3.0in]{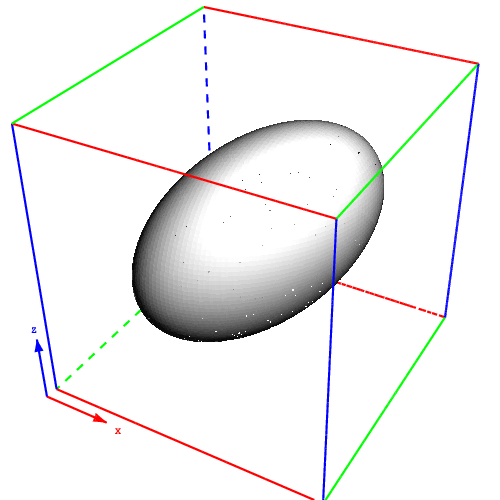} 
  \caption{The rotated $(1,.8..6)$ ellipsoid} 
  \label{fig:ex1}
\end{figure}

\begin{figure}
  \centering
  \includegraphics[width=3.0in]{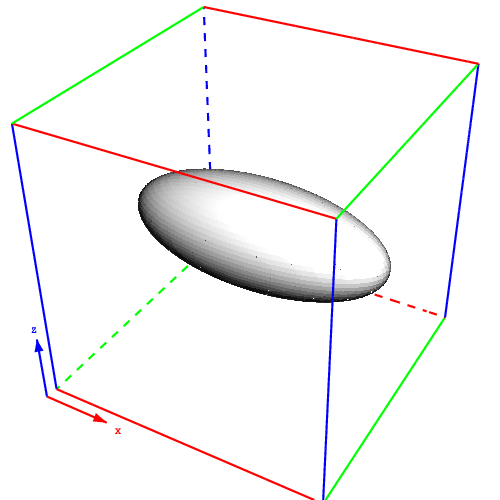} 
  \caption{The $(1,.4,.4)$ ellipsoid} 
  \label{fig:ex2}
\end{figure}

\begin{figure}
  \centering
  \includegraphics[width=3.0in]{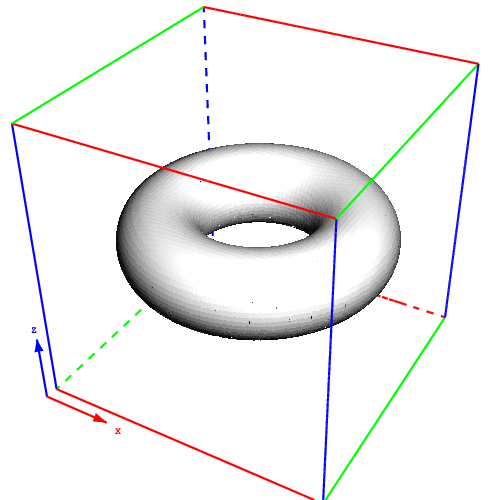} 
  \caption{The torus} 
  \label{fig:ex3}
\end{figure}

\begin{figure}
  \centering
  \includegraphics[width=3.0in]{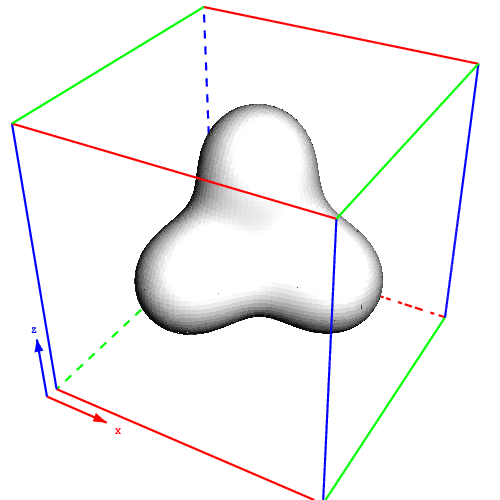} 
  \caption{The four-atom molecular surface} 
  \label{fig:ex4}
\end{figure}

\begin{figure}
  \centering
  \includegraphics[width=3.0in]{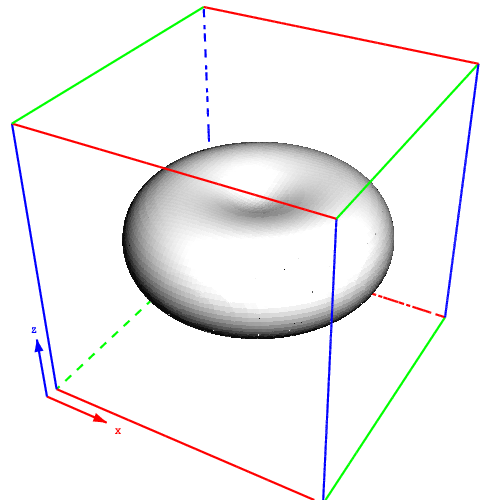} 
  \caption{The Cassini oval surface} 
  \label{fig:ex5}
\end{figure}

The errors for each example are presented in the tables with
$N = 64$, $128$ or $256$ and with $\delta/h = 1$, $2$, or $3$.  Both $L^2$
and maximum, or $L^\infty$, errors are given.  They are displayed first
for the irregular grid points, then for the regular grid points, and
finally for the quadrature nodes on the surface.  In the $L^2$ norms,
points are given equal weight, and thus for the quadrature nodes this
measure effectively gives extra weight to the overlap regions.  In all the
examples the angle $\theta$
in the partition of unity on the sphere is $70^o$.  The final table gives
the number of quadrature nodes on each surface with $N = 256$.

For the irregular points the smallest errors are generally for $\delta/h =1$,
while accuracy approaching $O(h^3)$ is observable with $\delta/h = 2$ or $3$.
(We have found that errors are larger with $\delta/h < 1$.)
The errors at the regular points is $O(h^2)$ as expected.  The errors at the
quadrature points are generally smaller than those at the irregular points for
$\delta/h = 2$ or $3$, but not for $\delta/h = 1$; this reflects
 the fact
that the method of Section 3.3
for evaluation on the surface improves the smoothing error directly but the discretization
error is improved indirectly by the smoothing.  With $\delta/h = 3$ the errors on the
surface decrease rapidly with refinement.  We repeated the computations on the surface
with angle $60^o$ and found similar but slightly larger errors.
To test further refinement, we computed $u$ on the surface for Example 1
with $N = 512$.  With angle $70^o$ and $\delta/h = 2$, the $L^2$ and
$L^\infty$ errors were 3.03E-8 and 7.43E-7; with $\delta/h = 3$ they were
3.09E-10 and 1.33E-8.  With angle $60^o$, they were
 6.83E-10, 3.82E-8 for $\delta/h = 2$ and 2.42E-10, 1.09E-8 for $\delta/h = 3$.
We conclude that for practical use we can reliably choose $\delta/h = 2$
for values at points near the surface and $\delta/h = 3$
with the special method of Sec. 3.3 for points on the surface.

In these calculations we used analytical values of $g^{ij}$ since the surface
was specified.  However these could easily be computed from grid values of
a level set function (see Appendix B) and are needed only to $O(h)$ accuracy.
We used analytical values of the densities $\varphi$ and $\psi$ at the quadrature
points and computed other values numerically from these.  We did so because
in solving an integral equation we would only know values at the quadrature points.
Occasionally the computation of $\partial^{(k)}_r\varphi$ in (\ref{N2}) failed for lack of
nearby quadrature points near the edge of the support of one $\zeta^{k,\theta}$.
In such a case we set the contribution to zero since $\zeta^{k,\theta}$
is very small there.  We also treated the ellipsoid of Example 1 without rotation;
the errors were similar to those displayed with the rotation but slightly smaller.
We tried an ellipsoid thinner than in Example 2; we found that the interpolation
stencil needed for the corrections $\T_1$, $\N_1$ failed with $N=64$ but worked with larger $N$.

\begin{table}[ht]

\caption{Errors for the rotated $(1,.8,.6)$ ellipsoid.
In each table the angle $\theta$ is $70^o$ and $h = 2.2/N$
for the $N^3$ grid.
$L^2$ and $L^\infty$ errors are displayed for (i) irregular grid points
near the surface, (ii) regular grid points, and (iii) quadrature points
on the surface.} 
\label{tab:tabrot}
\vspace{6pt}
\centering 
\begin{tabular}{| c | c | c | c | c | c | c | c |} \hline
$\delta$ & grid 
& $\|e_h^{\rsup{irreg}}\|_{2}$ & $\|e_h^{\rsup{irreg}}\|_{\infty}$ 
& $\|e_h^{\rsup{reg}}\|_{2}$ & $\|e_h^{\rsup{reg}}\|_{\infty}$ 
& $\|e_h^{\rsup{quad}}\|_{2}$ & $\|e_h^{\rsup{quad}}\|_{\infty}$
 \\ \hline  \hline
\multirow{3}{*}{$h$}
& $64^3$ & 3.32E-5 & 2.57E-4 & 3.29E-5 & 5.05E-4 & 9.80E-5 & 1.24E-3\\ \cline{2-8}
& $128^3$ & 4.14E-6 & 3.54E-5 & 7.96E-6 & 1.07E-4 & 2.10E-5 & 3.33E-4  \\ \cline{2-8}
& $256^3$ & 9.91E-7 & 6.55E-6 & 2.10E-6 & 2.92E-5 & 5.36E-6 & 8.69E-5  \\ \hline
\hline
\multirow{3}{*}{$2h$} 
& $64^3$ & 1.39E-4 & 8.64E-4 & 1.82E-4 & 2.39E-3 & 3.15E-5 & 4.32E-4  \\ \cline{2-8}
& $128^3$ & 1.78E-5 & 1.14E-4 & 4.35E-5 & 4.67E-4 & 1.71E-6 & 3.94E-5 \\ \cline{2-8}
& $256^3$ & 3.33E-6 & 1.39E-5 & 1.13E-5 & 1.36E-4 & 1.38E-7 & 4.69E-6   \\ \hline
\hline
\multirow{3}{*}{$3h$} 
& $64^3$ & 4.78E-4 & 2.91E-3 & 3.65E-4 & 5.06E-3 &  2.59E-5 & 2.83E-4 \\ \cline{2-8}
& $128^3$ & 5.75E-5 & 3.82E-4 & 8.79E-5 & 9.45E-4 & 1.13E-6 & 1.68E-5  \\ \cline{2-8}
& $256^3$ & 7.84E-6 & 4.78E-5 & 2.27E-5 & 2.70E-4 & 1.67E-8 & 6.36E-7 \\ \hline

\end{tabular}
\end{table}

\begin{table}[ht]

\caption{Errors for the $(1,.4,.4)$ ellipsoid} 
\label{tab:tabthin}
\vspace{6pt}
\centering 
\begin{tabular}{| c | c | c | c | c | c | c | c |} \hline
$\delta$ & grid 
& $\|e_h^{\rsup{irreg}}\|_{2}$ & $\|e_h^{\rsup{irreg}}\|_{\infty}$ 
& $\|e_h^{\rsup{reg}}\|_{2}$ & $\|e_h^{\rsup{reg}}\|_{\infty}$ 
& $\|e_h^{\rsup{quad}}\|_{2}$ & $\|e_h^{\rsup{quad}}\|_{\infty}$ \\
\hline\hline
\multirow{3}{*}{$h$}
 & $64^3$ & 4.46E-5 & 3.27E-4 & 2.62E-5 & 4.39E-4 & 1.48E-4 & 9.26E-4\\ \cline{2-8}
& $128^3$ & 8.53E-6 & 8.73E-5 & 6.79E-6 & 1.21E-4 &  3.25E-5 & 3.07E-4 \\ \cline{2-8}
& $256^3$ & 1.08E-6 & 1.39E-5 & 1.69E-6 & 3.15E-5 &  7.35E-6 & 6.70E-5\\ 
\hline\hline
\multirow{3}{*}{$2h$}
& $64^3$ & 3.29E-4 & 2.33E-3 & 1.69E-4 & 2.05E-3    & 4.95E-5 & 2.66E-4\\ \cline{2-8}
& $128^3$ & 4.29E-5 & 2.73E-4 & 4.07E-5 & 4.17E-4 & 6.82E-6 & 9.36E-5  \\ \cline{2-8}
& $256^3$ & 5.11E-6 & 3.57E-5 & 1.04E-5 & 1.05E-4 & 4.04E-7 & 1.00E-5 \\
\hline\hline
\multirow{3}{*}{$3h$}
& $64^3$ & 1.12E-3 & 6.85E-3 & 3.37E-4 & 4.29E-3    & 4.59E-5 & 2.35E-4 \\ \cline{2-8}
& $128^3$ & 1.43E-4 & 9.62E-4 & 8.07E-5 & 7.85E-4 & 5.22E-6 & 5.20E-5 \\ \cline{2-8}
& $256^3$ & 1.76E-5 & 1.23E-4 & 2.09E-5 & 2.22E-4 & 2.03E-7 & 3.91E-6 \\
\hline

\end{tabular}
\end{table}

\begin{table}[ht]

\caption{Errors for the torus} 
\label{tab:tabtorus}
\vspace{6pt}
\centering 
\begin{tabular}{| c | c | c | c | c | c | c | c |} \hline
$\delta$ & grid 
& $\|e_h^{\rsup{irreg}}\|_{2}$ & $\|e_h^{\rsup{irreg}}\|_{\infty}$ 
& $\|e_h^{\rsup{reg}}\|_{2}$ & $\|e_h^{\rsup{reg}}\|_{\infty}$ 
& $\|e_h^{\rsup{quad}}\|_{2}$ & $\|e_h^{\rsup{quad}}\|_{\infty}$
 \\ \hline  \hline
\multirow{3}{*}{$h$}
& $64^3$& 7.19E-5 & 3.57E-4 & 5.17E-5 & 5.03E-4 & 1.48E-4 & 1.29E-3  \\ \cline{2-8}
& $128^3$ & 8.61E-6 & 7.56E-5 & 9.04E-6 & 1.14E-4 & 3.14E-5 & 3.05E-4\\ \cline{2-8}
& $256^3$ & 1.02E-6 & 9.61E-6 & 2.16E-6 & 2.01E-5 & 7.02E-6 & 7.50E-5   \\ 
\hline\hline
\multirow{3}{*}{$2h$} 
& $64^3$ & 2.42E-4 & 7.94E-4 & 1.74E-4 & 1.42E-3 & 8.08E-5 & 4.16E-4  \\ \cline{2-8}
& $128^3$ & 2.89E-5 & 9.54E-5 & 4.41E-5 & 3.31E-4 & 8.53E-6 & 8.34E-5 \\ \cline{2-8}
& $256^3$ & 3.52E-6 & 1.25E-5 & 1.12E-5 & 8.85E-5 & 4.76E-7 & 7.46E-6 \\
\hline\hline
\multirow{3}{*}{$3h$} 
& $64^3$ & 8.17E-4 & 2.68E-3 & 3.17E-4 & 3.01E-3 & 6.35E-5 & 2.80E-4  \\ \cline{2-8}
& $128^3$ & 9.92E-5 & 3.28E-4 & 8.51E-5 & 6.93E-4 &  7.05E-6 & 4.85E-5\\ \cline{2-8}
& $256^3$ & 1.22E-5 & 4.22E-5 & 2.18E-5 & 1.84E-4 &  2.46E-7 & 2.80E-6\\ 
\hline

\end{tabular}
\end{table}

\begin{table}[ht]

\caption{Errors for the molecular surface} 
\label{tab:tabmole}
\vspace{6pt}
\centering 
\begin{tabular}{| c | c | c | c | c | c | c | c |} \hline
$\delta$ & grid 
& $\|e_h^{\rsup{irreg}}\|_{2}$ & $\|e_h^{\rsup{irreg}}\|_{\infty}$ 
& $\|e_h^{\rsup{reg}}\|_{2}$ & $\|e_h^{\rsup{reg}}\|_{\infty}$ 
& $\|e_h^{\rsup{quad}}\|_{2}$ & $\|e_h^{\rsup{quad}}\|_{\infty}$
 \\ \hline  \hline
\multirow{3}{*}{$h$}
& $64^3$ & 6.84E-5 & 4.16E-4 & 4.02E-5 & 4.81E-4 &  1.61E-4 & 1.47E-3 \\ \cline{2-8}
& $128^3$ & 5.98E-6 & 5.55E-5 & 1.09E-5 & 1.50E-4 & 3.24E-5 & 3.24E-4  \\ \cline{2-8}
& $256^3$ & 1.03E-6 & 1.30E-5 & 2.73E-6 & 3.96E-5 &  7.60E-6 & 8.33E-5 \\ 
\hline\hline
\multirow{3}{*}{$2h$} 
& $64^3$ & 3.30E-4 & 1.53E-3 & 2.51E-4 & 3.73E-3 &  6.80E-5 & 6.07E-4\\ \cline{2-8}
& $128^3$ & 3.99E-5 & 1.81E-4 & 6.02E-5 & 7.37E-4 &  3.43E-6 & 6.78E-5 \\ \cline{2-8}
& $256^3$ & 4.96E-6 & 2.36E-5 & 1.45E-5 & 1.57E-4 &  2.32E-7 & 4.34E-6 \\ 
\hline\hline
\multirow{3}{*}{$3h$} 
& $64^3$ & 1.11E-3 & 5.12E-3 & 5.08E-4 & 8.56E-3 &  6.35E-5 & 4.35E-4\\ \cline{2-8}
& $128^3$ & 1.38E-4 & 6.26E-4 & 1.18E-4 & 1.58E-3 &  2.01E-6 & 3.46E-5 \\ \cline{2-8}
& $256^3$ & 1.72E-5 & 8.17E-5 & 2.80E-5 & 3.27E-4 &   5.40E-8 & 1.24E-6\\ 
\hline

\end{tabular}
\end{table}

\begin{table}[ht]

\caption{Errors for the Cassini oval surface} 
\label{tab:taboval}
\vspace{6pt}
\centering 
\begin{tabular}{| c | c | c | c | c | c | c | c |} \hline
$\delta$ & grid 
& $\|e_h^{\rsup{irreg}}\|_{2}$ & $\|e_h^{\rsup{irreg}}\|_{\infty}$ 
& $\|e_h^{\rsup{reg}}\|_{2}$ & $\|e_h^{\rsup{reg}}\|_{\infty}$ 
& $\|e_h^{\rsup{quad}}\|_{2}$ & $\|e_h^{\rsup{quad}}\|_{\infty}$
 \\ \hline  \hline
\multirow{3}{*}{$h$}
& $64^3$ & 4.87E-5 & 2.94E-4 & 3.40E-5 & 3.60E-4 & 1.20E-4 & 1.01E-3  \\ \cline{2-8}
& $128^3$ & 3.78E-6 & 3.07E-5 & 7.25E-6 & 6.59E-5 & 2.37E-5 & 2.27E-4 \\ \cline{2-8}
& $256^3$ & 6.82E-7 & 5.75E-6 & 1.82E-6 & 1.68E-5 & 5.62E-6 & 5.75E-5\\
\hline\hline
\multirow{3}{*}{$2h$} 
& $64^3$ & 2.02E-4 & 8.64E-4 & 1.55E-4 & 1.20E-3 & 5.18E-5 & 3.53E-4 \\ \cline{2-8}
& $128^3$ & 2.46E-5 & 1.20E-4 & 3.80E-5 & 3.38E-4 & 3.15E-6 & 3.76E-5  \\ \cline{2-8}
& $256^3$ & 3.10E-6 & 1.56E-5 & 9.59E-6 & 8.17E-5 & 2.19E-7 & 3.63E-6\\
\hline\hline
\multirow{3}{*}{$3h$} 
& $64^3$ & 6.83E-4 & 2.65E-3 & 2.97E-4 & 2.58E-3 & 4.47E-5 & 2.20E-4 \\ \cline{2-8}
& $128^3$ & 8.61E-5 & 3.91E-4 & 7.35E-5 & 7.50E-4 &  1.84E-6 & 1.86E-5\\ \cline{2-8}
& $256^3$ & 1.08E-5 & 5.13E-5 & 1.86E-5 & 1.78E-4 & 4.65E-8 & 6.96E-7 \\
\hline

\end{tabular}
\end{table}

\begin{table}[ht]

\caption{The number of quadrature nodes with $N = 256$} 
\vspace{6pt}
\centering 
\begin{tabular}{| c | c | c | c | c | c |} \hline
Example & 1 & 2 & 3 & 4 & 5 \\ \hline
Number of nodes & 144388 & 70790 & 142168 & 126789 & 133014 \\  \hline

\end{tabular}
\end{table}

\section{\label{sec_discussion}Discussion}

The numerical results illustrate the performance of the method and are in general
agreement with the qualitative predictions.  They show that reasonable accuracy
can be obtained with moderate resolution, and the observed order of accuracy gives
confidence that the errors will reduce with further refinement.
Here we comment on possible improvements.  

As noted in Section 3, the discretization error depends on the angle $\theta$ in
the partition of unity on the unit sphere, defined in Sec. 2.  We need
$\theta > 55^o$ to cover the sphere.  As $\theta$ increases toward $90^o$ 
we expect the accuracy to deteriorate because of the dependence of the discretization error on $\theta$, as explained in Sec. 3.4. 
Here we used $\theta = 70^o$ as a compromise between the extremes.
In our experiments the errors were not very sensitive to the choice of angle for
$60^o \leq \theta \leq 80^o$.
We found slightly larger errors with $\theta = 60^o$ than for $70^o$.
A possible explanation is that the gradient of the partition of unity functions is
larger for the smaller angle.  It is unclear whether this can be improved or is
an inherent limitation.

In the discretization corrections $\T2$, $\N2$ we summed over
$n = (n_1,n_2)$ with $|n_j| \leq 20$, but the number of terms actually needed is much
smaller.  In fact for  $\delta/h \geq 2$ these corrections are usually negligible.  
They could be modified to include an estimate of the number of terms needed to avoid
unnecessary work. 

The treecode of \cite{Duan:2000} cannot be used directly with the method of
Section 3.3 for evaluation on the surface because of the differences in the regularized
kernels.  However, the treecode could be applied in this case by modifying the
recurrence formulas for the Taylor coefficients as derived in \cite{Duan:2000}.
This could be done in future work.

Other than \cite{Duan:2000}, fast summation methods have not been developed
specifically for regularized kernels.  
Among existing codes, one possible alternative 
that might be used in the present computations is the
kernel-independent fast multipole method (KIFMM) of 
L. Ying, G. Biros, D. Zorin \cite{LYing:2004}.
We have calculated examples using this code, even though it was not intended for
regularized kernels.  We found difficulty maintaining good accuracy, especially with
larger $\delta/h$, perhaps because the regularization degrades the accuracy of the
linear problems solved in the KIFMM.  We emphasize that this
is a use of the KIFMM for which it was not intended. A
summation method of fast multipole type designed particularly for these regularized
kernels could improve the efficiency of this method without loss of accuracy.

\section*{Acknowledgments}

Research of the first author was supported in part by the National Science 
Foundation of the USA under Grant DMS--1312654.
Research of the second author was supported in part by the National Natural 
Science Foundation of China under Grants DMS--11101278 and DMS--91130012.
Research of the second author was also supported in part by the Young 
Thousand Talents Program of China. 
Research of the third author was supported by the National Science 
Foundation of the USA under Grant DMS-0806482.

\appendix

\section{Regularization error for the single layer potential}

The regularization correction for the single layer potential,
evaluated at a point near the surface, can be derived using
the method of \cite{Beale:2004}, Sec. 2.  With $G_\del$ and $G$
as in (\ref{Gdelta}),
we approximate the error in the single layer potential with
density $\psi$, evaluated at a point $\vect x$ near the surface.
Since the error is local, we write it as an integral in
one coordinate patch,
regarding $\psi$ as a function of $\al = (\al_1,\al_2)$,
\be
 \eps = \int \left[ G_\del(\vect y(\al) - \vect x) - 
G(\vect y(\al)- \vect x)\right]
            \psi(\al)\,dS(\al) \,.
\ee
For simplicity, we will assume that $\vect x$ is along
the normal line from $0 \in \Gamma$, so that $\vect x = b \vect n_0$ for some $b$,
where $\vect n_0$ is the unit normal at $0$.  We also assume the coordinates
are chosen so that $\al(0) = 0$,  $g_{ij}(\al) =  \del_{ij} + O(|\al|^2)$,
and the tangent vectors $T_1$, $T_2$ have the directions of principal curvature.
Thus
\be
 \eps =  \frac{1}{4\pi}\int \frac{\text{erfc}(r/\del)}{r}\psi(\al)\,dS(\al) \,, \qquad
           r = |y(\al) - x|
\ee

Proceeding as before, we make a near-identity coordinate change $\al \to \xi$
such that $|\xi|^2 + b^2 = r^2$.  We get
 \be
   \eps =  \frac{1}{4\pi}\int \frac{\text{erfc}(\sqrt{|\xi|^2 + b^2}/\del)}
          {(|\xi|^2 + b^2)^{1/2}}  w(\xi,b)\,d\xi
 \ee
with
\be
   w(\xi,b) = \psi\left|\frac{\pa\al}{\pa\xi}\right||T_1\times T_2|
\ee
We will see that we can neglect terms in $w$ of the form $O(|\xi|^2 + b^2)$.
We can approximate $\psi$ in $\xi$, with leading term $\psi_0 = \psi(0)$,
\be
   \psi = \psi_0 + \psi_j(1 + bq/2)\xi_j + 
          \textstyle{\frac12} \psi_{ij}\xi_i\xi_j + O(|\xi|^3) + O(b^3) \,.
\ee
where $q = \kappa_1\xi_1^2/|\xi|^2 + \kappa_2\xi_2^2/|\xi|^2$
and $\kappa_1, \kappa_2$ are the principal curvatures at $0$.  
For the other two factors in $w$, we have
\be
   \det(\pa\al/\pa\xi) = 1 + bq + O(|\xi|^2) + O(b^2) \,.
\ee
and
\be
 |T_1\times T_2| = 1 + O(|\xi|^2) \,.
 \ee

In the $\xi$-integral for $\eps$, the odd part of $w$ will contribute zero.
Thus we can replace $w$ with an approximation to its even part.
Combining the three factors above, we get
\be
  w^{even}(\xi,b) =  \psi_0(1 + bq) + O(|\xi|^2 + b^2)
\ee
We now substitute in the integral, change to
polar coordinates, and substitute
$|\xi| = \del s$ and $b = \del\lam$ to obtain
\be \eps =  \frac{\del}{2}\psi_0(1 + \del\lam H)
           \int_0^\infty \frac{\text{erfc}(\sqrt{s^2 + \lam^2})}
          {(s^2 + \lam^2)^{1/2}} s\,ds + O(\del^3)
 \ee
with $H = (\kappa_1 + \kappa_2)/2$, the mean curvature.
With $r = \sqrt{s^2 + \lam^2}$, and $s\,ds = r\,dr$, the integral simplifies to
\be I = \int_0^\infty \frac{\text{erfc}(r)}
          {r} s\,ds =  \int_{|\lam|}^\infty \text{erfc}(r)\,dr
        = e^{-\lam^2}/\sqrt{\pi} - |\lam|\text{erfc}{|\lam|}
\ee
and finally
\be \eps = \frac{\del}{2}\psi_0(1 + \del\lam H)
                   \left(e^{-\lam^2}/\sqrt{\pi} - |\lam|\text{erfc}{|\lam|}   \right)
    + O(\del^3)
\ee
leading to the correction (\ref{T1}).


\section{\label{sec_monge}Formulas for Monge Patches}

We summarize formulas needed for the corrections of Sec. \ref{sec_layer_potentials}
when applied in a coordinate system such as $x_3 = f(x_1,x_2)$,
often called a Monge patch.
Given $f$, let $f_i = \partial f/\partial x_i$, $i = 1,2$, and similarly
for a second derivative $f_{ij}$.
The metric tensor $(g_{ij})$ and its inverse $(g^{ij})$ are 
$$  
(g_{ij}) = \left(  \begin{array}{cc}
    1+f_1^2 & f_1f_2 \\ f_1f_2 & 1+f_2^2  \end{array} \right)\,, \qquad
(g^{ij}) = \frac{1}{g} 
   \left(  \begin{array}{cc}
      1+f_2^2 & -f_1f_2 \\-f_1f_2 & 1+f_1^2  \end{array} \right) . 
$$
where $g = \det(g_{ij}) = 1 + f_1^2 + f_2^2$.
The Gauss curvature is
$$ K = g^{-2}\left(f_{11}f_{22} - f_{12}^2\right)$$
The mean curvature is
$$  
	H = \pm \tfrac12 g^{-3/2}\left[(1+f_2^2)f_{11} + (1+f_1^2)f_{22} 
	- 2f_1f_2f_{12}\right] 
$$
where the sign is $+$ if $x_3 > f(x_1,x_2)$ outside and $-$ otherwise.

The surface Laplacian has the general formula
$$ 
\laplace_S \varphi = \sum_{i,j} \frac{1}{\sqrt{g}} \, 
\partial_j\left(\sqrt{g} \, g^{ij} \partial_i \varphi \right)
   = \sum_{i,j} g^{ij} \left(\partial_i\partial_j \varphi \right) 
       + \sum_i c_i \partial_i \varphi   $$
where
$$          c_i = \sum_j \frac{1}{\sqrt{g}} \partial_j\left(\sqrt{g} g^{ij} \right) 
                =  \sum_j \frac{1}{\sqrt{g}} \partial_j\left(\frac{1}{\sqrt{g}} 
     \left( g g^{ij}\right) \right) $$
With some calculation we find
$$  \frac{1}{\sqrt{g}} \, \partial_j\left(\frac{1}{\sqrt{g}} \right) = 
               - g^{-2}\left(f_1f_{1j} + f_2f_{2j}\right)  $$
and subsequently
$$   c_i = 
           g^{-2}f_i\left[2f_1f_2f_{12} - (1+f_2^2)f_{11} - (1+f_1^2)f_{22} \right] 
           = \mp \frac{2}{\sqrt{g}} H f_i  $$

Suppose the surface is defined by $\phi(x_1,x_2,x_3) = 0$, with $\phi > 0$ outside.
Near a given point, there is at least one Monge patch; suppose
we can solve for $x_3 = f(x_1,x_2)$ as above.
By differentiating implicitly we get $f_i = - \phi_i/\phi_3$ etc.
We can use these to express $g_{ij}$ and  $g^{ij}$.
A more convenient expression for the mean curvature $H$
on a surface $\{\phi =  0\}$ is based on the classical formula
$$ 2H = - \nabla\cdot n =  - \nabla\cdot(\nabla \phi / |\nabla \phi|) $$
If we carry out the differentiation
we get 
$$ 2H = - |\nabla\phi|^{-3} \left( \phi_{ii}\phi_j^2 - \phi_i\phi_j\phi_{ij} \right) $$
summed over $i,j$.  After canceling and combining terms, we obtain
a formula such as in \cite{Osher:2003}, p. 12.


\bibliography{bywcicp.bbl}

\end{document}